\documentclass{article}
\usepackage[utf8]{inputenc}
\usepackage[left=1in,right=1in,bottom=1in]{geometry}
\usepackage[citestyle=alphabetic,
            bibstyle=alphabetic,
            maxbibnames=10,
            giveninits=true,
            doi=false,
            url=false,
            isbn=false,
            sorting=nyt,
            natbib=true]{biblatex}
\usepackage{amsmath}
\usepackage{amssymb}
\usepackage{mathtools}
\usepackage{amsthm}
\usepackage{dsfont}
\usepackage{subcaption}
\usepackage{xcolor}
\usepackage{setspace}
\usepackage{url}
\usepackage{cases}
\usepackage{tablefootnote}
\usepackage{authblk}

\newtheorem{assum}{Assumption}
\newtheorem{thm}{Theorem}

\newtheorem{cor}{Corollary}
\newtheorem{rem}{Remark}

\newcommand{\E}[2]{\mathbb{E}_{ #1 }\left[ #2 \right]}
\newcommand{\A}{\mathcal{A}}
\newcommand{\bracket}[2]{\left \langle #1 , #2 \right \rangle}
\DeclareMathOperator*{\argmin}{arg\,min}

\author[1,2]{Flemming Holtorf\thanks{E-mail: \texttt{holtorf@mit.edu}}}
\author[1,3]{Alan Edelman\thanks{E-mail: \texttt{edelman@mit.edu}}}
\author[1,4,5]{Christopher Rackauckas\thanks{E-mail: \texttt{crackauc@mit.edu}}}
\affil[1]{Computer Science and Artificial Intelligence Laboratory, Massachusetts Institute of Technology, 77 Massachusetts Ave, 02139 Cambridge, MA}
\affil[2]{Department of Chemical Engineering, Massachusetts Institute of Technology, 77 Massachusetts Ave, 02139 Cambridge, MA}
\affil[3]{Department of Mathematics, Massachusetts Institute of Technology, 77 Massachusetts Ave, 02139 Cambridge, MA}
\affil[4]{JuliaHub, Alewife, MA, USA}
\affil[5]{Pumas-AI, Centreville, VA, USA}

\title{Stochastic Optimal Control via Local Occupation Measures}

\begin{document}

\maketitle

\begin{abstract}                         
Viewing stochastic processes through the lens of occupation measures has proved to be a powerful angle of attack for the theoretical and computational analysis of stochastic optimal control problems. We present a simple modification of the traditional occupation measure framework derived from resolving the occupation measures locally on a partition of the control problem's space-time domain. This notion of local occupation measures provides fine-grained control over the construction of structured semidefinite programming relaxations for a rich class of stochastic optimal control problems with embedded diffusion and jump processes via the moment-sum-of-squares hierarchy. As such, it bridges the gap between discretization-based approximations to the Hamilton-Jacobi-Bellmann equations and occupation measure relaxations. We demonstrate with examples that this approach enables the computation of high quality bounds for the optimal value of a large class of stochastic optimal control problems with significant performance gains relative to the traditional occupation measure framework. 
\end{abstract}

\section{Introduction}
The optimal control of stochastic processes is one of the archetypical problems of decision-making under uncertainty with a myriad of applications in science and engineering. Despite their ubiquity, however, only a small subset of such stochastic optimal control problems admits the computation of a globally optimal control policy in a tractable and certifiable manner. As a consequence, engineers are often forced to resort to one of many available heuristics for the design of control policies in practice. And although such heuristics often perform remarkably well, they seldom come with a simple mechanism to quantify rigorously the degree of suboptimality they introduce, ultimately leaving it to the engineer's intuition when the controller design process shall be terminated.

In response to this undesirable situation, the task of computing theoretically guaranteed bounds gauging the best attainable performance for various classes of stochastic optimal control and related problems has received considerable attention in the recent past; contributions range from bounding schemes for the optimal control of systems governed by deterministic nonlinear ordinary \cite{lasserre2008nonlinear,henrion2008nonlinear,gaitsgory2009linear,henrion2023occupation} and partial differential equations \cite{korda2018nonlinear_pdes,korda2022moments,henrion2023infinite}, to discrete-time Markov control problems \cite{hernandez1999linear,savorgnan2009discrete}, to the control of diffusion and other continuous-time stochastic processes \cite{helmes2000numerical,cho2002linear,lamperski2018analysis,holtorf2023performance}. 
The underlying ideas have further found application in devising bounding schemes for the statistics of exit and stationary distributions of uncontrolled processes, for instance in mathematical finance \cite{lasserre2006pricing,henrion2023revisiting} and chemical physics \cite{ghusinga2017exact,dowdy2018bounds,backenkohler2019bounding,holtorf2024tighter}. In particular the framework of occupation measures has proved to be a versatile and effective approach for the construction of such bounding schemes. The notion of occupation measures allows for the translation of a rich class of stochastic optimal control and analysis problems into infinite-dimensional linear programs (LPs) over Borel measure spaces~\cite{fleming1989convex,vinter1993convex,bhatt1996occupation,kurtz1998existence}. Using Hausdorff moment conditions~\cite{helmes2001computing} or the moment-sum-of-squares (MSOS) hierarchy~\cite{lasserre2001global,parrilo2000structured}, these infinite-dimensional LPs in turn admit finite and tractable conic programming relaxations. In particular, the semidefinite programming (SDP) relaxations furnished by the MSOS hierarchy have been demonstrated to be computationally favorable and useful~\cite{lasserre2004sdp}. A key limitation of this framework, however, remains in its poor scalability. Specifically, the problem size of the SDP relaxations grows combinatorially with the hierarchy level while often high levels are necessary to establish informative bounds in practice. This challenge is further compounded by the notorious numerical ill-conditioning of Hankel moment matrices involving high-order momens~\cite{tyrtyshnikov1994bad,beckermann2000condition,papp2011optimization}.

In this paper, we propose a simple modification to the traditional occupation measure framework to address these practical limitations. Specifically, we consider localizing the occupation measures on a joint partition of the control problem's state space and time horizon. In contrast to previous works exploring similar notions of localized occupation measures for the optimal control of switched~\cite{claeys2016modal} and hybrid systems~\cite{abdalmoaty2013measures}, our key contribution lies in a refined analysis of the associated measure transport equations. By explicitly accounting for causality and continuity properties of the process paths, we establish transport equations that generate a richer set of moment constraints than those employed in the traditional framework and in turn furnish tighter conic relaxations via the MSOS hierarchy. Notably, these relaxations can be further tightened by refining the spatio-temporal partition rather than increasing the hierarchy level, offering two key practical advantages:
\begin{enumerate}
    \item evasion of numerical ill-conditioning due to consideration of higher-order moments considered when increasing the level in the MSOS hierarchy
    
    \item more fine-grained and interpretable control over tightening of the SDP relaxations when compared to increasing the level in the MSOS hierarchy  
\end{enumerate}
As we demonstrate with examples, these advantages hold the potential to construct equally or even tighter relaxations that can be solved notably faster than those derived with the traditional approach. Another potential advantage worth mentioning, though beyond the scope of this work, is that the proposed approach is similar in spirit to a wide range of discretization techniques for the numerical solution of partial differential equations; as such, the resultant MSOS relaxations exhibit a similar weakly-coupled block structure which may be exploited, for example with distributed optimization techniques~\cite{Shin2020DecentralizedProblems,Boyd2010DistributedMultipliers}. 

A partitioning approach closely related to the here proposed local occupation measure framework has recently been studied by Cibulka et al.~\cite{cibulka2021spatio} in the context of approximating the region of attraction for deterministic control systems via sum-of-squares programming. In another related work, Holtorf and Barton~\cite{holtorf2024tighter} have used temporal partitioning in order to improve MSOS bounding schemes for trajectories of stochastic chemical systems modeled by jump processes. Both works report significant computational merits of the respective modifications. Here, we unify and extend both works by introducing the notion of local occupation measures which applies beyond deterministic control problems to jump and diffusion control problems alike. The resultant framework is independent from and can be complemented by other approaches aimed at improving the tractability and practicality of the MSOS hierarchy, such as symmetry reduction ~\cite{riener2013exploiting,augier2023symmetry}, sparsity exploitation~\cite{schlosser2020sparse,wang2021exploiting,zheng2019sparse}, and linear/second-order cone programming hierarchies~\cite{ahmadi2014dsos,ahmadi2017optimization,ahmadi2019construction,ahmadi2019dsos}.

The remainder of this article is structured as follows: In Section \ref{sec:background}, we review the concept of occupation measures and show how it enables the construction of tractable convex relaxations for a large class of stochastic optimal control problems with embedded diffusion processes. In Sections \ref{sec:dual} and \ref{sec:primal}, we introduce the notion of local occupation measures and study its interpretation in the context of stochastic optimal control from the dual (polynomial) and primal (moment) perspective, respectively. Section \ref{sec:scaling} is dedicated to highlight the advantages of the proposed framework with regard to the scaling properties and structure of the resultant SDP relaxations. In Section \ref{sec:case_study}, we showcase the potential of the proposed approach with an example problem from population control. In Section \ref{sec:extensions}, we discuss the extension of the described local occupation measure framework to discounted infinite horizon control problems as well as the control of jump processes, supported with an example from systems biology. We conclude with some final remarks in Section \ref{sec:conclusion}.

\section{Problem description \& preliminaries}\label{sec:background}
We consider a continuous-time diffusion process $x_t$ in $\mathbb{R}^{n_x}$ driven by a standard $\mathbb{R}^m$-Brownian motion $B_t$ and controlled by a non-anticipative control process $u_t$ in $\mathbb{R}^{n_u}$,
\begin{align}
    {\rm d}x_t = f(x_t,u_t) \, {\rm d}t + g(x_t,u_t) \, {\rm d} B_t, \label{eq:SDE}
\end{align}
and study the associated finite horizon optimal control problem
\begin{align}
    J \coloneqq \inf_{u_t} \quad &\E{\nu_0}{\int_{[0,T]} \ell(x_t,u_t) \, {\rm d}t + \phi(x_T)} \tag{OCP} \label{eq:OCP}\\
    \text{s.t.} \quad & x_t \text{ satisfies \eqref{eq:SDE} on } [0,T] \text{ with } x_0 \sim \nu_0, \nonumber\\
    & (x_t, u_t) \in X \times U \text{ on } [0,T],\nonumber \\
    & u_t \text{ is non-anticipative.}\nonumber
\end{align}
Here, $\mathbb{E}_{\nu_0}$ denotes the expectation with respect to the probability measure $\mathbb{P}_{\nu_0}$ over the paths of the diffusion process \eqref{eq:SDE}. The subscript $\nu_0$ indicates the dependence on the distribution of the initial state, which we assume to be known. Throughout, we further assume that all problem data is described in terms of polynomials in the following sense.
\begin{assum}\label{ass:poly}
The drift coefficient $f:X \times U \to \mathbb{R}^{n_x}$, diffusion matrix $gg^\top:X \times U \to \mathbb{R}^{n_x \times n_x}$, stage cost $l:X\times U \to \mathbb{R}$ and terminal cost $\phi:X \times U \to \mathbb{R}$ are componentwise polynomial functions jointly in both arguments. The state space $X$ and the set of admissible control actions $U$ are basic closed semialgebraic sets. 
\end{assum}

We say a control process $u_t$ is admissible if the the controlled process $(x_t,u_t)$ satisfies the constraints in Problem \eqref{eq:OCP}. Furthermore, we make the following well-posedness assumption that ensures that the optimal value of \eqref{eq:OCP} is finite.
\begin{assum}\label{ass:finite_moments}
    The controlled diffusion process \eqref{eq:SDE} has finite moments for any admissible control process, i.e., $\E{\nu_0}{p(x_t, u_t)}$ is finite for all polynomials $p$ and $t \in [0,T]$.  
\end{assum}
Note that this assumption does not impose strong practical restrictions as it is for instance implied if the distribution of the controlled process has exponentially decaying tails or if $X$ and $U$ are compact. 

The key insight enabling the construction of convex relaxations of \eqref{eq:OCP} is that the controlled process described by \eqref{eq:SDE} admits a weak-form characterization in terms of a pair of occupations measures: the instantaneous and expected state-action occupation measure \cite{fleming1989convex,kurtz1998existence,bhatt1996occupation}. This characterization endows the control problem with a convex, albeit infinite-dimensional, geometry, sidestepping the nonlinear dependence of the paths of the diffusion process \eqref{eq:SDE} on the control process.

The instantaneous occupation measure $\nu$ is given by the probability to observe $x_T$ in any Borel set $B\subset X$. Formally, we define 
\begin{align*}
    \nu(B) \coloneqq \mathbb{P}_{\nu_0}\left[ x_T \in B \right].
\end{align*}
 or equivalently,
\begin{align*}
    \bracket{w}{\nu} \coloneqq \E{\nu_0}{w(T,x_T)}
\end{align*}
for every continuous test function $w \in \mathcal{C}([0,T]\times X)$, where 
\begin{align*}
    \bracket{w}{\nu} \coloneqq \int_{X} w(T,x) \, {\rm d}\nu(x)
\end{align*}
denotes the standard duality bracket between continuous functions and finite measures.

The expected state-action occupation measure $\xi$ is defined as the average time the controlled process $(t,x_t,u_t)$ remains in a Borel subset of $[0,T] \times X \times U$; formally, we define
\begin{align*}
    \xi(B_T \times B_X \times B_U)  \coloneqq \E{\nu_0}{\int_{[0,T] \cap B_T} \mathds{1}_{B_X\times B_U}((x_t,u_t)) \, {\rm d}t}
\end{align*}
for any Borel subsets $B_T \subset [0,T]$, $B_X \subset X$, $B_U \subset U$; or equivalently, 
\begin{align*}
    \bracket{w}{\xi} \coloneqq \E{\nu_0}{\int_{[0,T]} w(t,x_t,u_t)\, {\rm d}t }
\end{align*}
for any continuous test function $w \in \mathcal{C}([0,T]\times X\times U)$. The instantaneous and expected state-action occupation measures are finite, non-negative measures by construction. 

The occupation measure pair $(\nu, \xi)$ characterizes the expected time evolution of sufficiently smooth observables $w \in \mathcal{C}^{1,2}([0,T]\times X)$ of the process\footnote{that is functions on the domain $[0,T]\times X$ with continuous first and second derivatives (in the sense of Whitney~\cite{whitney1992analytic}) in the first and second argument, respectively.} by Dynkin's formula \cite[Theorem 1.24]{oksendal2007applied}, 
\begin{align*}
    \E{\nu_0}{w(T, x_T)} = \E{\nu_0}{w(0, x_0)} + \E{\nu_0}{\int_{[0,T]} \mathcal{A} w(s,x_s,u_s)\,{\rm d}s},
\end{align*}
or equivalently,
\begin{align}
    \bracket{w}{\nu} = \bracket{w}{\nu_0} + \bracket{\A w}{\xi}, \label{eq:mom_Dynkin}
\end{align}
where $\A:\mathcal{C}^{1,2}([0,T] \times X) \to \mathcal{C}([0,T] \times X \times U)$ denotes the (extended) infinitesimal generator of the diffusion process \eqref{eq:SDE} \cite{oksendal2007applied}, i.e.,
\begin{align*}
    \A:w(t,x) \mapsto \frac{\partial w}{\partial t}(t,x) + f(x,u)^\top \nabla_x w(t,x) + \frac{1}{2} \text{Tr}\left(gg^\top(x,u)  \nabla_x^2 w(t,x)\right).
\end{align*}
Conversely, we say that a measure pair $(\nu, \xi)$ is a weak solution to \eqref{eq:SDE} on the interval $[0,T]$ if it satisfies Equation \eqref{eq:mom_Dynkin} for all test functions $w \in \mathcal{C}^{1,2}([0,T] \times X)$. This notion of weak solutions to \eqref{eq:SDE} motivates the following weak form of \eqref{eq:OCP} \cite{fleming1989convex}:
\begin{align}
    J^* \coloneqq \inf_{\nu, \xi} \quad &\bracket{\ell}{\xi} + \bracket{\phi}{\nu} \tag{weak OCP} \label{eq:weakOCP}\\
    \text{s.t.} \quad & \bracket{w}{\nu} - \bracket{\A w}{\xi} =  \bracket{w}{\nu_0}, \ \forall w \in \mathcal{C}^{1,2}([0,T] \times X), \nonumber \\
    &\nu \in \mathcal{M}_+(X), \nonumber \\
    &\xi \in \mathcal{M}_+([0,T]\times X \times U). \nonumber
\end{align}
where $\mathcal{M}_+(Y)$ denotes the cone of finite, positive Borel measures supported on the set $Y$. Problem \eqref{eq:weakOCP} is an infinite-dimensional LP \cite{nash1987linear} and generally a strict relaxation of \eqref{eq:OCP}; in particular, $J > J^*$ may hold. 
Absence of a gap between $J^*$ and $J$ can however be guaranteed under suitable regularity conditions. Bhatt and Borkar~\cite{bhatt1996occupation} establish conditions for equivalence between \eqref{eq:OCP} and \eqref{eq:weakOCP} in the framework of relaxed (measure-valued) controls. 
Further, equivalence between the relaxed and strict control problem holds under additional convexity assumptions on the set of admissible controls and the dependence of drift coefficient, diffusion matrix, and stage cost on the control action~\cite{kushner1975existence,karoui1987compactification,haussmann1990existence}; for more recent results on related questions about occupation measure relaxations of deterministic optimal control and variational problems, see \cite{korda2022gap,henrion2024occupation}.

From a practical perspective, \eqref{eq:weakOCP} remains intractable as an infinite-dimensional LP; however, Assumption \ref{ass:poly} enables the construction of a sequence of increasingly tight SDP relaxations via the MSOS hierarchy \cite{parrilo2000structured,lasserre2001global}. To that end, \eqref{eq:weakOCP} is relaxed to the optimization over moment sequences of the measures $\nu$ and $\xi$ truncated at finite order $d$. For polynomial test functions, constraints of the form \eqref{eq:mom_Dynkin} reduce to affine constraints on the moment sequence as $\mathcal{A}$ maps polynomials to polynomials under Assumption \ref{ass:poly}. Similarly, the conic constraints $\nu \in \mathcal{M}_+(X)$ and $\xi \in \mathcal{M}_+([0,T]\times X \times U)$ can be relaxed to positive semidefiniteness constraints of certain moment and localizing matrices, which under Assumption \ref{ass:poly} reduce to linear matrix inequalities \cite{lasserre2001global}.

The infinite-dimensional LP dual \cite{nash1987linear} to \eqref{eq:weakOCP} has an informative interpretation that serves as motivation for the partitioning strategy presented in the next section. The dual reads 
\begin{align}
    \sup_{w} \quad & \int_{X} w(0,x) \, {\rm d}\nu_0(x) \tag{subHJB}\label{eq:subHJB}\\
    \text{s.t.} \quad & \A w + \ell \geq 0, \text{ on } [0,T] \times X \times U, \label{eq:dynamics_1} \\
    & w(T,\cdot) \leq  \phi, \text{ on } X, \label{eq:transversality_1} \\
    & w \in \mathcal{C}^{1,2}([0,T] \times X), \nonumber
\end{align}
where the decision variable $w$ can be interpreted as a smooth underestimator of the value function associated with the control problem \eqref{eq:OCP}:
\begin{cor}\label{cor:subsolution}
    Let $w$ be feasible for \eqref{eq:subHJB} and let $\delta_z$ denote the Dirac measure centered at $z$. Then, $w$ underestimates the value function
    \begin{align}
        V(t,z) \coloneqq \inf_{u_s} \ & \E{\delta_z}{\int_{t}^T \ell(x_s,u_s) \, {\rm d}s + \phi(x_T)} \label{eq:valuefunction}\\
                 \textup{s.t.} \ & x_s \text{ satisfies } \eqref{eq:SDE} \text{ on } [t, T] \text{ with } x_t \sim \delta_z, \nonumber \\
                    & (x_s, u_s) \in X\times U \text{ on } [t,T],\nonumber \\
                    & u_s \text{ is non-anticipative}. \nonumber
    \end{align}
    for any $(t,z)\in[0,T] \times X$.
\end{cor}
\begin{proof}
    Let $z\in X$ and $0 \leq t \leq T$ and fix any admissible control policy $u_s$, i.e., a control policy such that the path of the stochastic process $(x_s,u_s)$ remains in $X \times U$ on $[t,T]$. Then, Constraints \eqref{eq:dynamics_1} and \eqref{eq:transversality_1} imply that
    \begin{align*}
        \E{\delta_z}{-\int_{t}^T \A w(s,x_s,u_s) \, {\rm d}s + w(T,x_T)}  \leq  \E{\delta_z}{\int_{t}^T \ell(x_s,u_s) \, {\rm d}s + \phi(x_T)}.
    \end{align*}
    The left-hand-side coincides with $w(t,z)$ by Dynkin's formula. The result follows by minimizing over all admissible control policies. 
\end{proof}

Analogous to its primal counterpart, the MSOS hierarchy gives rise to a sequence of increasingly tight SDP restrictions of \eqref{eq:subHJB} by restricting $w$ to be a polynomial of degree at most $d$ and imposing the non-negativity constraints by means of sufficient sum-of-squares conditions \cite{lasserre2001global,parrilo2000structured}. The restriction is weakened by increasing the degree $d$ yielding a monotonically increasing sequence of lower bounds for the optimal value $J^*$ of \eqref{eq:weakOCP}. The following theorem establishes a set of easily verifiable conditions under which the limit point of this sequence is $J^*$ (implying also strong duality between \eqref{eq:subHJB} and \eqref{eq:weakOCP}).

\begin{thm}\label{thm:convergence}
    Let $J_d$ be the optimal value of the $d$\textsuperscript{th} level MSOS restriction of \eqref{eq:subHJB} (resp. relaxation of \eqref{eq:weakOCP}). If Assumption \ref{ass:poly} holds and moreover $X$ and $U$ are represented as
    \begin{align*}
        &X = \lbrace x : p_i(x) \geq 0, \ i = 1,\dots, v,   \ R_X - \| x\|^2_2 \geq 0 \rbrace, \\
        &U = \lbrace u : q_i(x) \geq 0, \ i = 1,\dots, w, \ R_U - \|u\|^2_2 \geq 0\rbrace,
    \end{align*}
    with suitable polynomials $p_i$ and $q_i$, and sufficiently large $R_X$ and $R_U$, then $J_d \uparrow J^*$.
\end{thm}
\begin{proof}
    First note that under the given assumptions, the set $[0,T]\times X\times U$ is compact. Thus, it suffices to impose condition \eqref{eq:mom_Dynkin} for all polynomial test functions in \eqref{eq:weakOCP} as it is a dense subset of $\mathcal{C}^{1,2}([0,T]\times X)$. Further observe that  constraint \eqref{eq:mom_Dynkin} implies that every feasible pair $(\nu, \xi)$ has constant mass. Specifically, for test functions $w(t,x)\equiv 1$ and $w(t,x)=t$, constraint \eqref{eq:mom_Dynkin} reduces to $\langle 1, \nu\rangle = 1$ and $\langle 1, \xi \rangle = T$, respectively. The result thus follows from \cite[Corollary 8]{tacchi2022convergence}.
\end{proof}
Finally, we emphasize that, while the MSOS bounds remain valid even for problems with unbounded state and action spaces, compactness of $X$ and $U$ as per the hypotheses of Theorem \ref{thm:convergence} is indeed necessary to guarantee absence of a gap between the limit point of the sequence of MSOS bounds and the best attainable control performance; see~\cite[Example 3.3]{lasserre2006pricing} for an example where an unbounded state space leads to a finite gap.

\section{The dual perspective revisited: piecewise polynomial approximation}\label{sec:dual}
In order to construct improved approximations to the value function in the spirit of \eqref{eq:subHJB}, we consider a generalization of problem \eqref{eq:subHJB} that seeks a {\em piecewise} smooth underapproximation of the value function over the problem's space-time domain $[0,T]\times X$. To that end, we consider a discretization $0 = t_0 < t_1 < \dots < t_{n_T}=T$ of the control horizon and a collection of state space restrictions $X_1,\dots,X_{n_X} \subset X$ which satisfy the following assumption and hence form a partition of $X$.
\begin{assum}\label{ass:partition}
    The collection $X_1,\dots, X_n \subset \mathbb{R}^{n_X}$ satisfies
    \begin{enumerate}
        \item $X = \cup_{k=1}^{n_X} X_k$,
        \item $X_i \cap X_j = \emptyset$ for all $1\leq i\neq j \leq n_X$. 
        \item the closure $\bar{X}_i$ and shared boundaries $\partial X_{ij} = \partial X_i \cap \partial X_j$ are basic closed semialgebraic for all $1 \leq i \neq j \leq 1,\dots,n_X$. 
    \end{enumerate} 
\end{assum} 
The elements $[t_{i-1}, t_i] \times X_k$ then form a partition of the problem's entire space-time domain and we can formulate the following natural generalization of \eqref{eq:subHJB}:
\begin{align}
    %
    %
    \sup_{w} \quad & \sum_{k = 1}^{n_X} \int_{X_k} w_{1,k}(0,x) \, {\rm d}\nu_0(x) \tag{pw-subHJB} \label{eq:disc_subHJB} \\
    \text{s.t.} \ \ \ & \A w_{i,k} + \ell \geq 0 \text{ on } [t_{i-1}, t_i] \times X_k \times U, \ \forall (i,k) \in P, \label{eq:path} \\[0.5em]
                            & w_{i,k}(t_{i-1}, \cdot) \geq w_{i-1,k}(t_{i-1},\cdot) \text{ on } X_k, \ \forall (i,k) \in P^{\circ}, \label{eq:time} \\[0.5em]
                            & w_{i,k} = w_{i,j} \text{ on } [t_{i-1}, t_i] \times (\partial X_j \cap  \partial X_k), \ \forall (i,j,k) \in \partial P, \label{eq:boundary} \\[0.5em]
                            & w_{n_T,k}(T,\cdot) \leq \phi \text{ on } X_k, \ \forall k \in \lbrace{1,\dots, n_X}\rbrace, \label{eq:transversality} \\[0.5em]
                            & w_{i,k} \in \mathcal{C}^{1,2}([0,T] \times X_k), \ \forall (i,k) \in P, 
\end{align}
with the index sets
\begin{align*}
    &P \coloneqq \left\lbrace (i,k) : 1 \leq i \leq n_T, 1 \leq k \leq n_X\right\rbrace, \\
    &P^{\circ} \coloneqq \lbrace (i,k) : 2 \leq i \leq n_T, 1\leq k \leq n_X\rbrace, \\
    &\partial P \coloneqq \lbrace (i,j,k): 1\leq i \leq n_T, 1 \leq k \neq j \leq n_X\rbrace.
\end{align*} 

The constraints in Problem \eqref{eq:disc_subHJB} ensure that a valid underestimator of the value function can be constructed from the function pieces $\lbrace w_{i,k}: (i,k) \in P\rbrace$ for all elements of the partition. As such, Problem \eqref{eq:disc_subHJB} yields a lower bound for the optimal value of \eqref{eq:OCP}. This is formalized in the following Corollary.  
\begin{cor}\label{cor:pw_subsolution}
    Let $\{w_{i,k}: (i,k) \in P \}$ be feasible for \eqref{eq:disc_subHJB} and define
\begin{align*}
    w(t,x) = w_{i,k}(t,x) \text{ with } i = \max \{ j : t \in [t_{j-1}, t_j] \} \text{ and } k \text{ such that } x \in X_k. 
\end{align*}
    Then, $w$ underestimates the value function $V$ as defined in Equation \eqref{eq:valuefunction}.
\end{cor}
\begin{proof}[Sketch]
    The idea is to split the paths of the process $(t,x_t,u_t)$ up into pieces during which it is confined to a single subdomain $[t_{i-1}, t_i] \times X_k \times U$. For each of those pieces an analogous argument as in Corollary \ref{cor:subsolution} applies to show that $w_{i,k}$ underestimates the value function for a process confined to the partition element $[t_{i-1}, t_{i}] \times X_k \times U$. Additionally, Constraints \eqref{eq:time} and \eqref{eq:boundary} ensure conservatism when the process crosses between different time intervals and subdomains of the state space, respectively. Specifically, Constraint \eqref{eq:time} enforces that $w(t,x_t)$ can at most decrease when traced backward in time across the boundary between the intervals $[t_{i}, t_{i+1}]$ and $[t_{i-1},t_{i}]$, ensuring that $w$ cannot cross $V$ at such time points. Similarly, Constraint \eqref{eq:boundary} imposes spatial continuity thus enforces that $w$ cannot cross $V$ when the process crosses spatial boundaries between partition elements. The formal argument is given in Appendix \ref{app:cor1}.
\end{proof}
\begin{rem}
From the detailed argument presented in Appendix \ref{app:cor1}, it is clear that piecewise smooth function underapproximators to the value function can be constructed by imposing alongside constraints \eqref{eq:path} and \eqref{eq:transversality} the requirement that the underapproximator increases at the boundaries between the spatio-temporal partition elements in expectation under the dynamics of the controlled process. As the process naturally crosses the temporal boundaries only in the direction of positive time, it is therefore sufficient to impose the monotonicity condition \eqref{eq:time}. In contrast, stochastic vibrations may lead to crossing of spatial boundaries between partition elements in any direction and thus the stronger continuity condition \eqref{eq:boundary} must be enforced. In the case of a deterministic process ($g \equiv 0$), this condition can however be further relaxed in a similar spirit as for the temporal boundaries. Cibulka et al. \cite{cibulka2021spatio} show in a similar context that in this case it suffices to impose 
\begin{align*}
    (w_{i,k}(t,x) - w_{i,j}(t,x)) n_{j,k} (x)^\top f(x,u) \geq 0, \ \forall (t,x,u) \in [t_{i-1},t_i] \times \partial X_{jk} \times U,
\end{align*}
where $n_{j,k} (x)$ denotes the normal vector of the boundary between $X_j$ and $X_k$ pointing from $X_j$ to $X_k$ at $x$. Intuitively, this condition enforces the monotonicity conditions $w_{i,k}(t,x) - w_{i,j}(t,x) \geq 0$ ($w_{i,k}(t,x) - w_{i,j}(t,x) \leq 0$) if the dynamics allow for crossing of the boundary $\partial X_{jk}$ at $x$ only in the direction from $X_j$ to $X_k$ ($X_k$ to $X_j$).
\end{rem}

\section{The primal perspective revisited: local occupation measures}\label{sec:primal}
In this section, we discuss the primal counterpart of the construction presented in the previous section. The primal counterpart of \eqref{eq:disc_subHJB} reads
\begin{align}
     \inf_{\nu, \xi, \pi} & \quad  \sum_{(i,k) \in P} \bracket{\ell}{\xi_{i,k}}  + \sum_{k=1}^{n_X} \bracket{\phi}{\nu_{n_T,k}} \tag{pw-weak OCP} \label{eq:disc_weakOCP}\\ 
    \text{s.t.} & \quad \bracket{w}{\nu_{i, k}} - \bracket{w}{\nu_{i-1,k}} = \bracket{\A w}{\xi_{i,k}} + \sum_{j \neq k} \bracket{w}{\pi_{i,j,k}}, \ \forall w \in \mathcal{C}^{1,2}([t_{i-1},t_i] \times X_k), \ \forall (i,k) \in P, \nonumber \\[0.5em]
     &\quad \nu_{i,k} \in \mathcal{M}_+(X_k), \quad \forall (i,k) \in P, \nonumber \\[0.5em]
     &\quad \xi_{i,k} \in \mathcal{M}_+([t_{i-1},t_i]\times X_k \times U), \quad \forall (i,k) \in P,\nonumber \\[0.5em]
     & \quad \pi_{i,j,k} = -\pi_{i,k,j} \in \mathcal{M}([t_{i-1},t_i] \times \partial X_{jk}), \ \forall (i,j,k) \in \partial P, \nonumber 
\end{align}
where $\mathcal{M}(Y)$ refers to the space of signed measures supported on $Y$. The decision variables in \eqref{eq:disc_weakOCP} can be interpreted as localized generalization of the occupation measure pair introduced in Section \ref{sec:background}. Specifically, the restriction of the expected state-action occupation measures $\xi$ to a subdomain $[t_{i-1}, t_i] \times X_k \times U$ from the partition generates the local state-action occupation measure $\xi_{i,k}$:
\begin{align*}
    \xi_{i,k}(B_T \times B_X \times B_U) = \xi((B_T \cap [t_{i-1}, t_i]) \times (B_X \cap X_k) \times B_U).
\end{align*}
 
Likewise, the local instantaneous occupation measures with respect to different time points $t_i$ and subdomains $X_k$ are given by the restriction of the instantaneous occupation measure at time $t_i$ to $X_k$, i.e.,
\begin{align*}
    \nu_{i,k}(B) = \mathbb{P}_{\nu_0}(x_{t_i} \in B \cap X_k).
\end{align*}
The measure $\pi_{i,j,k}$ in \eqref{eq:disc_weakOCP} takes the role of a slack variable and accounts for transitions of the process between the spatial subdomains $X_j$ and $X_k$ in the time interval $[t_{i-1}, t_i]$. Formally, $\pi_{i,j,k}$ can be defined by  
\begin{align*}
    \bracket{w}{\pi_{i,j,k}} \coloneqq  \E{\nu_0}{\sum_{n=1}^{N^{jk}_{+}} w\left(\tau_{n+}^{jk}, x_{\tau_{n+}^{jk}}\right) - \sum_{n=1}^{N^{jk}_{-}} w\left(\tau_{n-}^{jk}, x_{\tau_{n-}^{jk}}\right) },
\end{align*}
where $\tau_{n+}^{jk}$ and $\tau_{n-}^{jk}$ denote the $n$\textsuperscript{th} time points in $(t_{i-1}, t_i)$ at which the process transitions from subdomain $X_j$ into $X_k$ and vice versa, respectively. With these interpretations, we can observe that the equality constraints in \eqref{eq:disc_weakOCP} reduce to Dynkin's formula applied between the stopping times of leaving and entering a given subdomain $X_k$ in the time interval $[t_{i-1},t_i]$ (see Appendix \ref{app:slack} for a more detailed derivation). 

Finally, it is important to emphasize here that the above interpretation of the decision variables in \eqref{eq:disc_weakOCP} as local occupation measures shows that every feasible point for \eqref{eq:disc_weakOCP} generates a feasible point for \eqref{eq:weakOCP} with equal objective value via the assignment $\xi = \sum_{(i,k) \in P} \xi_{i,k}$ and $\nu = \sum_{k =1}^{n_X} \nu_{n_T,k}$. Analogously, any smooth function $w$ that is feasible for \eqref{eq:subHJB} generates upon restriction to the individual subdomains of the partition a feasible point for \eqref{eq:disc_subHJB} with equal objective value. 
This property carries over directly to the respective MSOS restrictions and relaxations, provided that every subdomain closure $\bar{X}_k$ is represented in terms of a strictly greater set of polynomial inequalities than $X$ is. Under this mild condition, the bounds furnished by the MSOS restrictions and relaxations of \eqref{eq:disc_subHJB} and \eqref{eq:disc_weakOCP} are therefore at least as tight, and often tighter, than those obtained from MSOS restrictions and relaxations of their traditional counterparts, even when the optimal values of the underlying infinite-dimensional LPs coincide.

\section{Moment-sum-of-squares approximations: structure \& scaling} \label{sec:scaling}
The construction of tractable relaxations of the problems \eqref{eq:subHJB} or \eqref{eq:weakOCP} relies on the restriction to optimization over polynomials of fixed degree $d$ or the relaxation to optimization over moment sequences truncated at order $d$, respectively. Increasing this approximation order $d$ has traditionally been the only mechanism used to weaken the restriction, respectively strengthen the relaxation, to improve the resultant bounds to a desired level. The main motivation behind the proposed partitioning approach lies in circumventing the limited practicality and interpretability of this tightening mechanism. With the proposed notion of local occupation measures, refinement of the space-time domain partition serves as an additional bound tightening mechanism. Table \ref{tab:scaling} summarizes how the MSOS SDP restrictions and relaxations of \eqref{eq:disc_subHJB} and \eqref{eq:disc_weakOCP} scale in size with respect to the different tightening mechanisms of increasing $n_X, n_T$ (refining the partition), or $d$ (increasing the approximation order). The linear scaling of the SDP sizes with respect to $n_X$ and $n_T$ underlines the fine-grained control over the tightening process via refinement of the partition. In particular, it opens the door to exploit problem specific insights such as the knowledge of critical parts of the (extended) state space $[0,T]\times X$ to be resolved more finely than others, to construct tighter relaxations without incurring a combinatorial increase in the number of partition elements. This flexibility and interpretability is in stark contrast to tightening the bounds by increasing the approximation order $d$ as translating such insights into specific moments to be constrained or polynomial basis elements to be considered for the value function approximator is significantly less straightforward. 
It is further worth emphasizing that not only the linear scaling with respect to $n_T$ and $n_X$ is desirable but in particular that the invariance of the linear matrix inequality (LMI) dimension promotes practicality due to the unfavorable scaling of interior point algorithms whose running time scales worse than cubically with respect to this quantity~\cite{jiang2023faster}.

\begin{table}[h]
    \caption{Scaling of MSOS SDP approximations with respect to different tightening mechanisms.\tablefootnote{here, $n = 1+ n_x + \max \lbrace \text{deg}_x f - 1, \text{deg}_x g - 2\rbrace$} }\label{tab:scaling}
    \centering
    \begin{tabular}{l|lll}
         & \#\textbf{variables} & \#\textbf{LMIs} & \textbf{LMI dimension}  \\
         \hline
         $d$ & $O({n + d \choose d})$ & $O(1)$ & $O({n + \left \lfloor d/2 \right \rfloor \choose \left \lfloor d/2 \right \rfloor})$ \\
         $n_T$ & $O(n_T)$ & $O(n_T)$ & $O(1)$\\
         $n_X$ & $O(n_X)$ & $O(n_X)$ & $O(1)$\\
         \hline
    \end{tabular}
\end{table}

Additionally, the problems \eqref{eq:disc_subHJB} and \eqref{eq:disc_weakOCP} give rise to highly structured SDPs. Specifically, all constraints involve only variables corresponding to adjacent subdomains. As a consequence, the structure of the constraints is analogous to those arising from discretized PDEs may be exploited with suitable distributed optimization algorithms and computing architectures. 
 
\section{Example: population control}\label{sec:case_study}

\subsection{Control problem}
We demonstrate the computational merits of the proposed local occupation measure framework with an example problem from the field of population control. The problem is adjusted from Savorgnan et al. \cite{savorgnan2009discrete} where it has been studied in a discrete time, infinite horizon setting. The objective is to control the population size of a primary predator and its prey in an ecosystem featuring the prey species, primary predator species as well as a secondary predator species. The interactions between the primary predator and prey population are described by a standard Lotka-Volterra model, while the effect of the secondary predator species is modeled by a Brownian motion. The population sizes are assumed to be controlled via hunting of the primary predator species. This model gives rise to the diffusion process
\begin{align*}
    &{\rm d}x_{t,1} =  (\gamma_1  x_{t,1} - \gamma_2  x_{t,1}  x_{t,2}) \, {\rm d}t + \gamma_5 x_{t,1} \, {\rm d} B_t,\\
    &{\rm d}x_{t,2} = (\gamma_4  x_{t,1}  x_{t,2} - \gamma_3  x_{t,2} - x_{t,2} u_t) \, {\rm d}t,
\end{align*}
where $x_1$, $x_2$, and $u$ refer to the prey species, predator species and hunting effort, respectively. The model parameters $\gamma = (1, 2, 1, 2, 0.025)$ and initial condition $x_0\sim \delta_{(1, 0.25)}$ are assumed to be known deterministically. Moreover, we assume that the admissible hunting effort is confined to $U = [0,1]$. Under these assumptions, it is easily verified that the process state $x_t$ evolves by construction within the non-negative orthant $X = \{x : x_1, x_2 \geq 0 \}$ for any admissible control policy. For the control problem we further choose a time horizon of $T = 10$ and stage cost
\begin{align*}
    \ell(x,u) = (x_1 - 0.75)^2 + \frac{(x_2 - 0.5)^2}{10} + \frac{(u - 0.5)^2 }{10} 
\end{align*}
penalizing variations from the target population sizes. 
\subsection{Partition of problem domain}
In order to investigate the effect of different discretizations on bound quality and computational cost, we utilize a simple grid partition of the state space $X$ as parameterized by the number of grid cells $n_1$ and $n_2$ in the $x_1$ and $x_2$ direction, respectively. As $X$ is the non-negative orthant in our example, and hence semi-infinite, we choose to discretize the compact interval box $[0,1.5]\times[0,1.5]$ with a uniform grid of $(n_1-1) \times (n_2-1)$ cells and cover the remainder of $X$ with appropriately chosen semi-infinite interval boxes. This choice is motivated by the insight that the uncontrolled system resides with high probability in $[0,1.5]\times[0.1.5]$.

The temporal domain is partitioned uniformly into $n_T$ subintervals, i.e., $t_i = i \Delta t$ with $\Delta t = T/n_T$. Throughout, we refer to a specific partition with the associated triple $(n_1, n_2, n_T)$. The computational experiments are conducted for all partitions corresponding to the triples $\{(n_1, n_2, n_T): 1 \leq n_1, n_2 \leq 5, 1\leq n_T \leq 10\}$.

\subsection{Evaluation of bound quality}\label{sec:boundquality}
In order to assess the tightness of the bounds obtained with different approximation orders and discretizations, we compare the relative optimality gap $(\bar{J}-\underbar{$J$}) / \bar{J}$, where $\underline{J}$ and $\bar{J}$ refer to the lower bound furnished by an instance of the MSOS restriction of \eqref{eq:disc_subHJB} and to the control cost associated with the {\em best known} admissible control policy, respectively. The best known control policy was constructed from the approximate value function $w^*$ obtained as the solution of the MSOS restriction of \eqref{eq:disc_subHJB} with approximation order $d = 4$ on the grid described by $n_1 = n_2 = 4$ and $n_T=10$. To that end, we employed the following control law mimicking a one-step model-predictive controller 
\begin{align*}
    u^*_{t} \in \argmin_{u \in U} \mathcal{A} w^*(t, x_t, u) + \ell(x_t,u)
\end{align*}
and estimated the associated control cost
\begin{align*}
    \bar{J} = \E{\nu_0}{\int_{[0,T]} \ell(x_t,u^*_{t}) \, {\rm d}t }
\end{align*}
by the ensemble average of $50,000$ sample trajectories.

\subsection{Computational aspects}
All computational experiments presented in this section were conducted on a MacBook M1 Pro with 16GB unified memory. All sum-of-squares programs and the associated SDPs were constructed using our custom developed and publicly available package \texttt{MarkovBounds.jl}\footnote{see \url{https://github.com/FHoltorf/MarkovBounds.jl}} built on top of \texttt{SumOfSquares.jl}~\cite{weisser2019polynomial} and the \texttt{MathOptInterface}~\cite{legat2022mathoptinterface}. All resultant SDPs were solved using Mosek v10.

\subsection{Results} 
\begin{figure}[!h]
    \centering
    \includegraphics[width=0.5\textwidth]{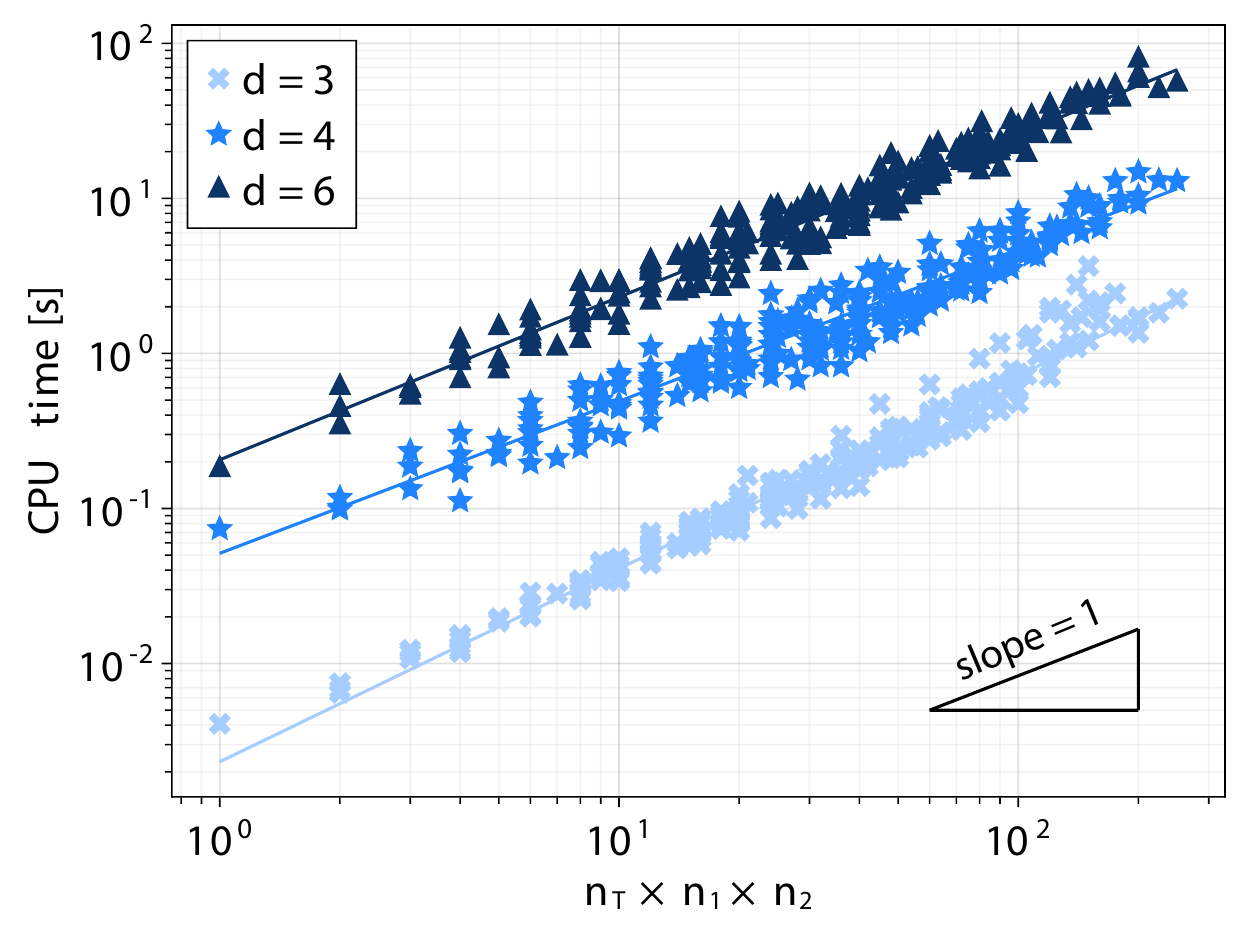}
    \caption{Linear scaling with respect to the the number of grid cells for fixed approximation order}
    \label{fig:linear_scaling}
\end{figure}

We put special emphasis on investigating the effect of refining the discretization of the problem domain on bound quality and computational cost. Focusing on the effect on computational cost in isolation first, Figure \ref{fig:linear_scaling} indicates that the computational cost for the solution of MSOS programs generated by the restriction of \eqref{eq:disc_subHJB} to polynomials of degree at most $d$ scales approximately linearly with the number of cells $n_1 \times n_2 \times n_T$ of the spatiotemporal partition. On the other hand, Figure \ref{fig:linear_scaling} also shows that increasing the approximation order $d$ results in a much more rapid increase in computational cost. These results are in line with the discussion in Section \ref{sec:scaling}.

\begin{figure}[!h]
    \centering
    \begin{subfigure}[t]{0.49\textwidth}
        \includegraphics[width=1.0\linewidth]{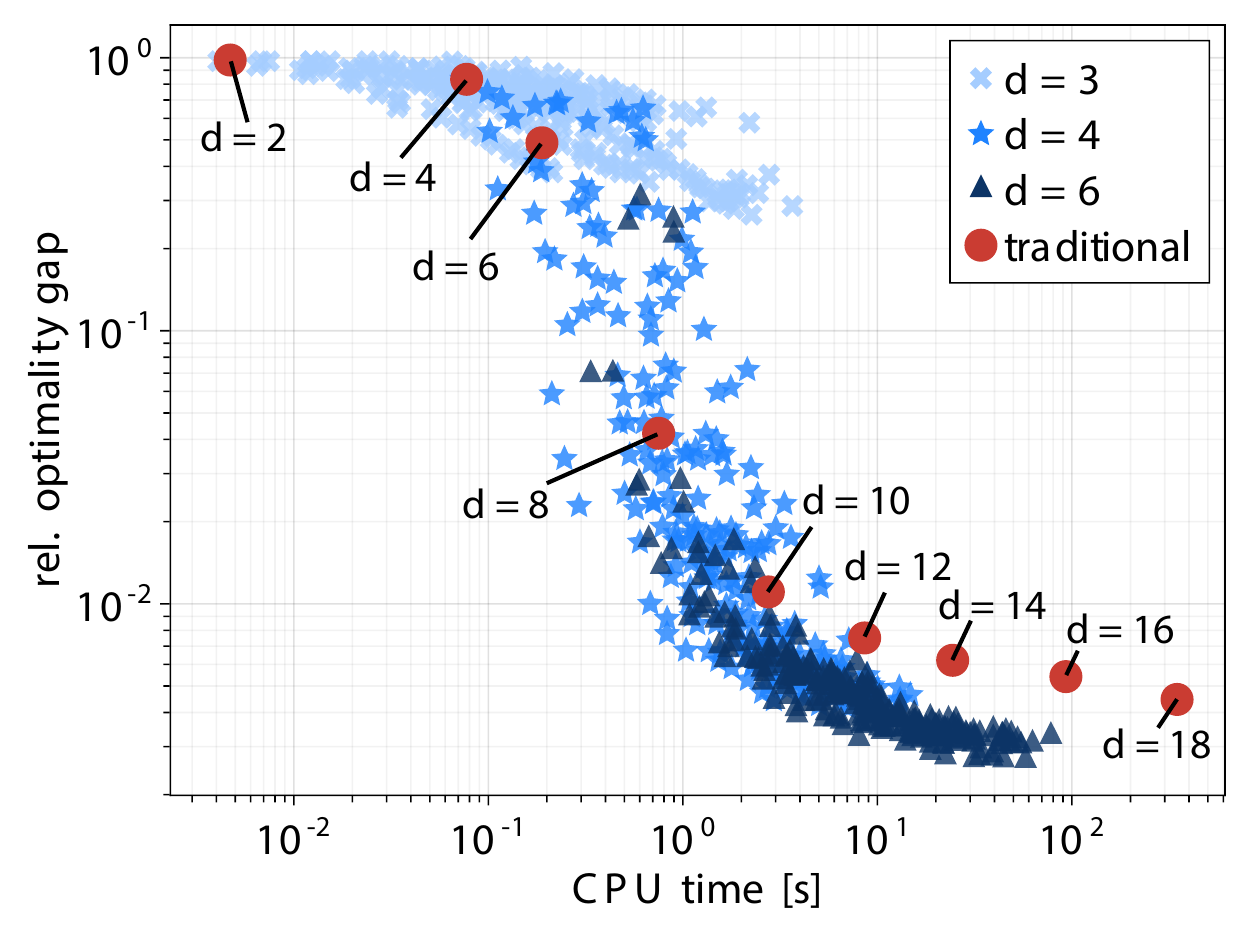}
    \caption{spatial \& temporal partitioning}
    \end{subfigure}
    \begin{subfigure}[t]{0.49\textwidth}
      \includegraphics[width=1.0\linewidth]{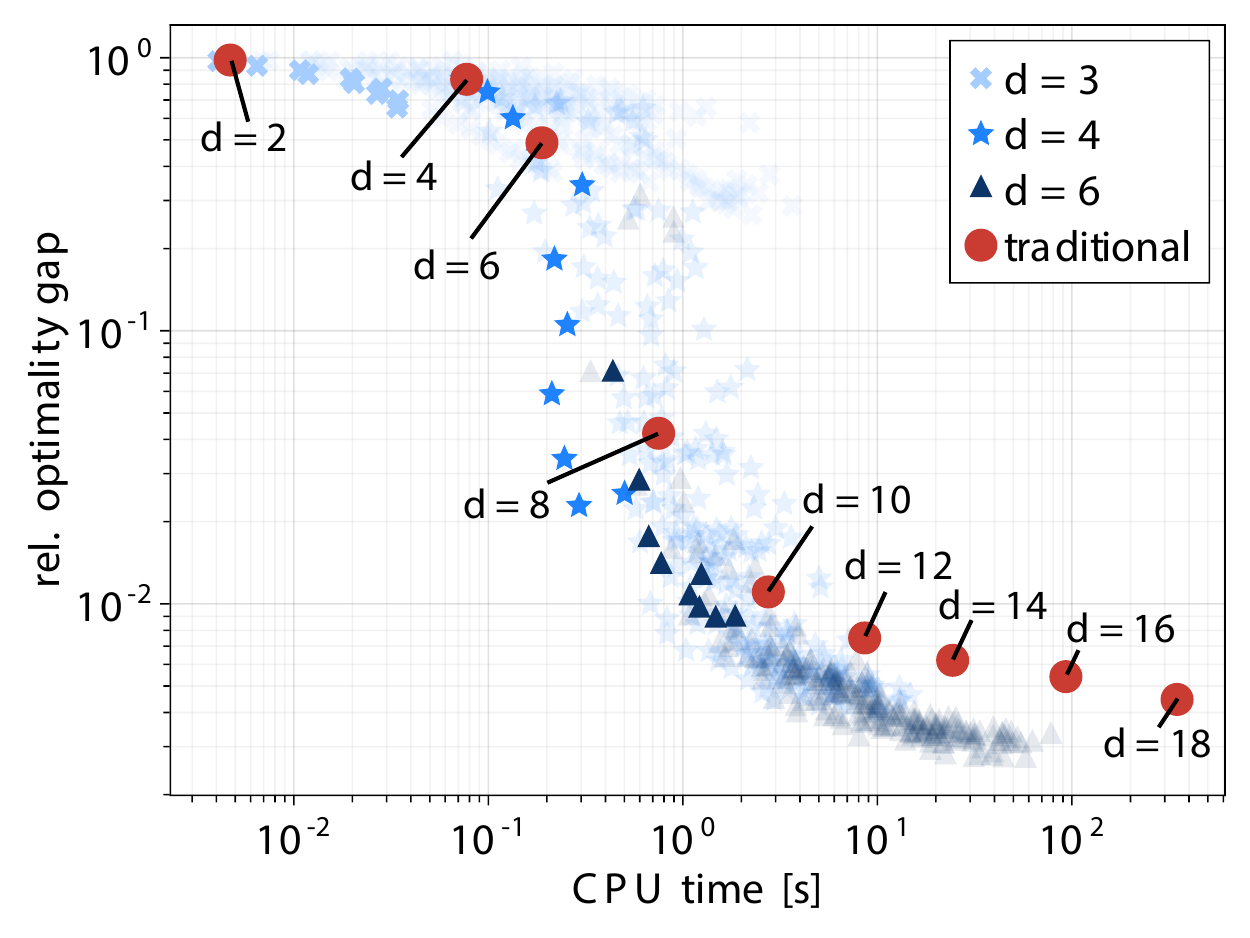}
      \caption{exclusively temporal partitions highlighted ($n_1 = n_2 = 1$)} \label{fig:temp}
    \end{subfigure}
    \caption{Trade-off between computational cost and bound quality for different approximation orders $d$ and domain discretizations $(n_1,n_2,n_T)$. The red markers correspond to MSOS restrictions of the labeled approximation order for the traditional formulation \eqref{eq:subHJB}.}
    \label{fig:performance}
\end{figure}

Figure \ref{fig:performance} shows the trade-off between bound quality and computational cost for different approximation orders and partitions. First, it is worth noting that the proposed partitioning strategy enables the computation of overall tighter bounds with an approximation order of only up to $d=6$ when compared to the traditional formulation with an approximation order of up to $d=18$. It is further worth emphasizing that beyond $d=18$, numerical issues prohibited an accurate solution of the SDPs arising from the traditional formulation such that no tighter bounds could be obtained this way. Furthermore, almost across almost the entire accuracy range a significant speed-up could be achieved by using the proposed partitioning strategy instead of only increasing the approximation order. Lastly, the results indicate that a careful choice of partitioning is crucial to achieve good performance. Figure \ref{fig:temp} suggests that for this example particularly good performance is achieved when only the time domain is partitioned; additionally partitioning the spatial domain becomes an effective means of bound tightening only after the time domain has been resolved sufficiently finely.

\section{Extensions}\label{sec:extensions}
Before we close, we briefly discuss two direct extensions to the described local occupation measure framework showcasing its versatility. 
\subsection{Discounted infinite horizon problems}
Consider the following discounted infinite horizon stochastic control problem with discount factor $\rho > 0$:
\begin{align*}
    \inf_{u_t} \quad &\E{\nu_0}{\int_{[0,\infty)} e^{-\rho t} \ell(x_t,u_t) \, {\rm d}t }\\
    \text{s.t.} \quad & x_t \text{ satisfies \eqref{eq:SDE} on } [0,\infty) \text{ with } x_0 \sim \nu_0, \\ 
    & (x_t, u_t) \in X \times U, \text{ on } [0,\infty), \\
    & u_t \text{ is non-anticipative}.
\end{align*}
The construction of a weak formulation of this problem akin \eqref{eq:weakOCP} can be done in full analogy to Section \ref{sec:background}. To that end, note that the infinitesimal generator $\A$ maps functions of the form $\hat{w}(t,x) = e^{-\rho t} w(t,x)$ to functions of the same form, i.e.,
\begin{align*}
    \A \hat{w}(t,x,u) = e^{-\rho t} (\A w(t,x,u) - \rho w(t,x,u)).
\end{align*}
By analogous arguments as in Section \ref{sec:background}, it therefore follows that any function $w \in C^{1,2}([0,\infty)\times X) $ that satisfies 
\begin{align*}
    \A w(t,x,u) - \rho w(t,x,u) + \ell(x,u) \geq 0, 
    \forall (t,x,u) \in [0,\infty) \times X \times U 
\end{align*}
generates a valid subsolution $\hat{w}(t,x) = e^{-\rho t} w(t,x)$ of the value function associated with the infinite horizon problem. Since the proposed partitioning approach does neither rely on boundedness of the state space nor control horizon in order to establish valid bounds, it follows that it readily extends to the infinite horizon setting.

\subsection{Stopped control problems}
So far we have treated $X$ as the intrinsic state space of the controlled process. When $X$ is to include path constraints that may be violated with non-zero probability under optimal control, problem \eqref{eq:OCP} must be extended to consider stopping the process upon leaving $X$:
\begin{align*}
    \inf_{u_t} \quad &\E{\nu_0}{\int_{[0,\tau]} \ell (x_t, u_t) \, dt + \phi(\tau, x_\tau)}\\
    \text{s.t.} \quad &x_t \text{ satisfies on } [0,\tau] \text{ with } x_0 \sim \nu_0, \\
    & \tau = \inf \lbrace t \in [0,T] : x_t \notin X \rbrace, \\
    & u_t \in U, \text{ on } [0,T]\\
    & u_t \text{ is non-anticipative}. \\
\end{align*}
By analogous arguments as in Section \ref{sec:background} any function $w \in \mathcal{C}^{1,2}([0,T]\times X)$ that satisfies
\begin{align*}
    &\mathcal{A}w + \ell \geq 0 \text{ on } [0,T] \times X \times U,\\
    &w \leq \phi \text{ on } \lbrace T \rbrace \times X \cup [0,T) \times \partial X
\end{align*}
underestimates the value function associated with the stopped problem. It follows that a primal-dual pair of infinite-dimensional LPs akin \eqref{eq:subHJB} and \eqref{eq:weakOCP} bounds the optimal value of the stopped problem from below. These LPs further admit tractable MSOS restrictions and relaxations when the terminal cost $\phi$ admits a representation
\begin{align*}
    \phi(t,x) = \begin{cases}
        \phi_{\partial X}(t, x), & (t,x) \in [0,T) \in \partial X\\
        \phi_{T}(x), &  (t,x) \in \lbrace T \rbrace \times X
    \end{cases}
\end{align*}
with polynomial pieces $\phi_{\partial X}$ and $\phi_T$. The local occupation measure framework may be applied analogous to Sections \ref{sec:dual} and \ref{sec:primal} to provide refinement.


\subsection{Jump processes with discrete state space}
Many application areas ranging from chemical physics to queuing theory call for models that describe stochastic transitions between discrete states. In those cases, jump processes are a common modeling choice \cite{gillespie1992rigorous,breuer2003markov}. In the following, we will show that the proposed local occupation measure framework extends with only minor modifications to stochastic optimal control of a large class of such jump processes. To that end, we will consider controlled, continuous-time jump processes driven by $m$ independent Poisson counters $n_{i}(t)$ with associated propensities $a_i(x_t,u_t)$:
\begin{align}
    {\rm d}x_t = \sum_{i=1}^{m} \left[ h_i(x_t, u_t) - x_t \right] \, dn_{i,t}. \label{eq:JDE}
\end{align}
We will again assume that the process can be fully characterized by polynomials, but we now additionally impose the assumption that the state space of the process is discrete.
\begin{assum}\label{ass:poly_jump}
    The jumps $h_i:X \times U \to X$, propensities $a_i:X \times U \to \mathbb{R}_+$, stage cost $l:X\times U \to \mathbb{R}$ and terminal cost $\phi:X \times U \to \mathbb{R}$ are polynomial functions jointly in both arguments. The state space is a discrete, countable set and the set of admissible control actions $U$ is a basic closed semialgebraic set.
\end{assum}
The local occupation measure framework outlined previously for diffusion processes can be extended for computing lower bounds on stochastic optimal control problems with such jump processes embedded:
\begin{align}
    \inf_{u_t} \quad &\E{\nu_0}{\int_{[0,T]} \ell(x_t,u_t) \, {\rm d}t + \phi(x_T)} \tag{jump OCP} \label{eq:jumpOCP} \\
    \text{s.t.} \quad & x_t \text{ satisfies \eqref{eq:JDE} on } [0,T] \text{ with } x_0 \sim \nu_0, \nonumber\\
    & (x_t, u_t) \in X \times U, \text{ on } [0,T],\nonumber\\
    & u_t \text{ is not anticipative}.\nonumber
\end{align}
Given the extended infinitesimal generator $\mathcal{A}: \mathcal{C}^{1,0}([0,T] \times X) \to \mathcal{C}([0,T]\times X \times U)$ associated with the process \eqref{eq:JDE},
\begin{align*}
    \mathcal{A} w \mapsto \frac{\partial w}{\partial t}(t,x) + \sum_{i=1}^m a_i(x,u) \left( w(t, h_i(x,u)) - w(t,x)\right),
\end{align*}
the weak form of \eqref{eq:jumpOCP} and its dual are analogous to \eqref{eq:weakOCP} and \eqref{eq:subHJB}, respectively. Further, given a partition of the problem's space-time domain as introduced in Section \ref{sec:dual}, the analog of Problem \eqref{eq:disc_subHJB} seeking a piecewise smooth subsolution of the value function takes the form
\begin{align}
     \sup_{w_{i,k} : (i,k) \in P} \quad & \sum_{k = 1}^{n_X} \int_{X_k} w_{1,k}(0,\cdot) \, {\rm d}\nu_0  \label{eq:discHJB} \tag{jump pw-subHJB} \\
    \text{s.t.} \ \ \qquad  & \A w_{i,k} + \ell \geq 0 \text{ on } [t_{i-1}, t_i] \times X_k \times U, \ \forall (i,k) \in P, \nonumber \\[0.5em]
                            & w_{i,k}(t_{i-1}, \cdot) \geq w_{i-1,k}(t_{i-1},\cdot) \text{ on } X_k, \ \forall (i,k) \in P^{\circ}, \nonumber \\[0.5em]
                            & w_{i,k} = w_{i,j} \text{ on } [t_{i-1}, t_i] \times N_{X_k}(X_j), \ \forall (i,j,k) \in \partial P,  \nonumber\\[0.5em]
                            & w_{n_T,k}(T,\cdot) \leq \phi \text{ on } X_k, \ \forall k \in \lbrace{1,\dots, n_X}\rbrace \nonumber \\[0.5em]
                            & w_{i,k} \in \mathcal{C}^{1,0}([0,T] \times X_k), \ \forall (i,k) \in P, \nonumber          
\end{align}
where $N_{X_k}(X_j)$ denotes the ``neighborhood'' of $X_k$ in $X_j$ defined as all states in $X_j$ which have a non-zero transition probability into $X_k$; formally,
\begin{align*}
    N_{X_k}(X_j) = \{ x \in X_j : \exists u \in U \text{ such that }  h_i(x, u) \in X_k \text{ for some } i \text{ and } a_i(x, u) > 0\}.
\end{align*}

Note that under Assumption \ref{ass:poly_jump}, $\mathcal{A}$ again maps polynomials to polynomials laying the basis for the application of the MSOS hierarchy to construct tractable relaxations of \eqref{eq:jumpOCP}. In contrast to the discussion in Section \ref{sec:background}, however, the state space $X$ of a jump process is closed basic semialgebraic if and only if it is finite. Thus, the MSOS hierarchy provides finite SDP relaxations of the weak form of \eqref{eq:jumpOCP} only in the case of a finite state space $X$. Moreover, even if $X$ is finite but of large cardinality, these relaxations may not be practically tractable due to the large number or high degree of the polynomial inequalities needed to describe such a set. If $X$ is infinite (or of sufficiently large cardinality), tractable MSOS relaxations can only be constructed at the price of introducing additional conservatism. From the dual perspective, this additional conservatism is introduced by imposing the non-negativity conditions in \eqref{eq:subHJB} on a basic semialgebraic overapproximation of $X$; in particular polyhedral overapproximations are a common choice \cite{dowdy2018bounds,ghusinga2017exact,kuntz2019bounding,sakurai2017convex,holtorf2024tighter}. 
The framework of local occupation measures provides a way to reduce this conservatism. 
In order to construct tractable relaxations for  \eqref{eq:discHJB} via the MSOS hierarchy, it of course remains still necessary to replace any infinite (or very large) partition element by a closed basic semialgebraic overapproximation; however, the union of suitably chosen overapproximations of the individual partition elements will generally be less conservative than a global overapproximation. 



\subsection{Example: optimal gene regulation for protein expression}
We demonstrate the efficacy of the local occupation measure framework for the control of jump processes with an example from cellular biology. Specifically, we consider the problem of optimal regulation of protein expression through actuation of the promoter kinetics in the biocircuit. The biocircuit is modeled as a jump process reflecting the stochastic nature of chemical reactions in cellular environments with low molecular copy numbers \cite{gillespie1992rigorous}. The associated jump process has three states encoding the molecular counts of protein ($x_1$), active promoter ($x_2$), and inactive promoter ($x_3$) undergoing jump transitions in response to the following chemical reactions with associated rates:
\begin{align}
    &\begin{array}{l} h_1 : (x_1, x_2, x_3) \mapsto (x_1 + 1, x_2, x_3), \\
       \hspace*{2.5cm} a_{1}(x, u) = 10 x_2 \end{array} \tag{expression} \\
    &\begin{array}{l} h_2 : (x_1, x_2, x_3) \mapsto (x_1 -1, x_2, x_3), \\
     \hspace*{2.5cm} a_{2}(x, u) = 0.1 x_1 \end{array} \tag{degradation} \\
    &\begin{array}{l} h_3 : (x_1, x_2, x_3) \mapsto (x_1, x_2-1, x_3+1), \\
     \hspace*{2.5cm} a_{3}(x, u) = 0.1 x_1 x_2 \end{array} \tag{repression} \\
    &\begin{array}{l} h_4 : (x_1, x_2, x_3) \mapsto (x_1, x_2 + 1, x_3 - 1), \\ 
     \hspace*{2.5cm} a_{4}(x,u) = 10 (1-u) x_3 \end{array} \tag{activation} \\
    &\begin{array}{l} h_5 : (x_1, x_2, x_3) \mapsto (x_1, x_2 - 1, x_3 + 1), \\ 
     \hspace*{2.5cm} a_5(x, u) = 10 u x_2 \end{array} \tag{inactivation} 
\end{align}
The expression of protein can be controlled indirectly via the activation and inactivation rates of the promoter. Admissible control actions $u$ are constrained to lie within the interval $U = [0,1]$. Moreover, we assume a deterministic initial condition $x_0 \sim \delta_{(0,1,0)}$ and exploit that due to the reaction invariant $x_{t,2} + x_{t,3} = x_{0,2} + x_{0,3}$ the state space $X$ is effectively two-dimensional, i.e., we eliminate $x_{t,3} = 1 - x_{t,2}$. It can be easily verified that, after elimination of the reaction invariant, the state space of the jump process is given by
\begin{align*}
   X = \{x \in \mathbb{Z}_+^2 : x_2 \in \{0, 1\} \}
\end{align*}
such that Assumption \ref{ass:poly_jump} is satisfied. 

The goal of the control problem is to stabilize the protein level in the cell at a desired value of 10 molecules. To that end, we choose to minimize the stage cost
\begin{align*}
    \ell(x,u) = (x_1 - 10)^2 + 10 (u - 0.5)^2
\end{align*}
over the horizon $[0,10]$. 

In order to investigate the effect of different partitions of the problem domain on bound quality and computational cost, we discretize the time horizon uniformly into $n_T$ intervals and partition the state space into $2n_X$ singletons 
\begin{align*}
    X_i = \begin{cases}
    \{(i-1, 0), & i \leq n_X\\
    \{(i-n_X - 1, 1)\}, & i > n_X
    \end{cases} \text{ for } i = 1,\dots, 2n_X
\end{align*}
and lump the remaining part of the state space in the last partition element
\begin{align*}
    X_{2n_X+1} = \{ x \in \mathbb{Z}_+^2: x_1 \geq n_X, x_2 \in \{0, 1\} \}.
\end{align*}
We explore the partitions corresponding to all combinations of $n_T \in \{2, 4, \dots, 18, 20\}$ and $n_X \in \{0, 8, \dots, 32, 40\}$.

Note that the partition elements $X_1, \dots, X_{2n_X}$ are already basic closed semialgebraic such that no overapproximation is required for the construction of valid MSOS restriction of the non-negativity constraints in \eqref{eq:discHJB}. In contrast, $X_{2n_X+1}$ is infinite, hence not basic closed semialgebraic. We therefore strengthen the formulation of the MSOS restriction of \eqref{eq:discHJB} by imposing the non-negativity conditions on the polyhedral convex hull of $X_{2n_X+1}$, thereby recovering tractability.

\begin{figure}
    \centering
    \includegraphics[width=0.5\textwidth]{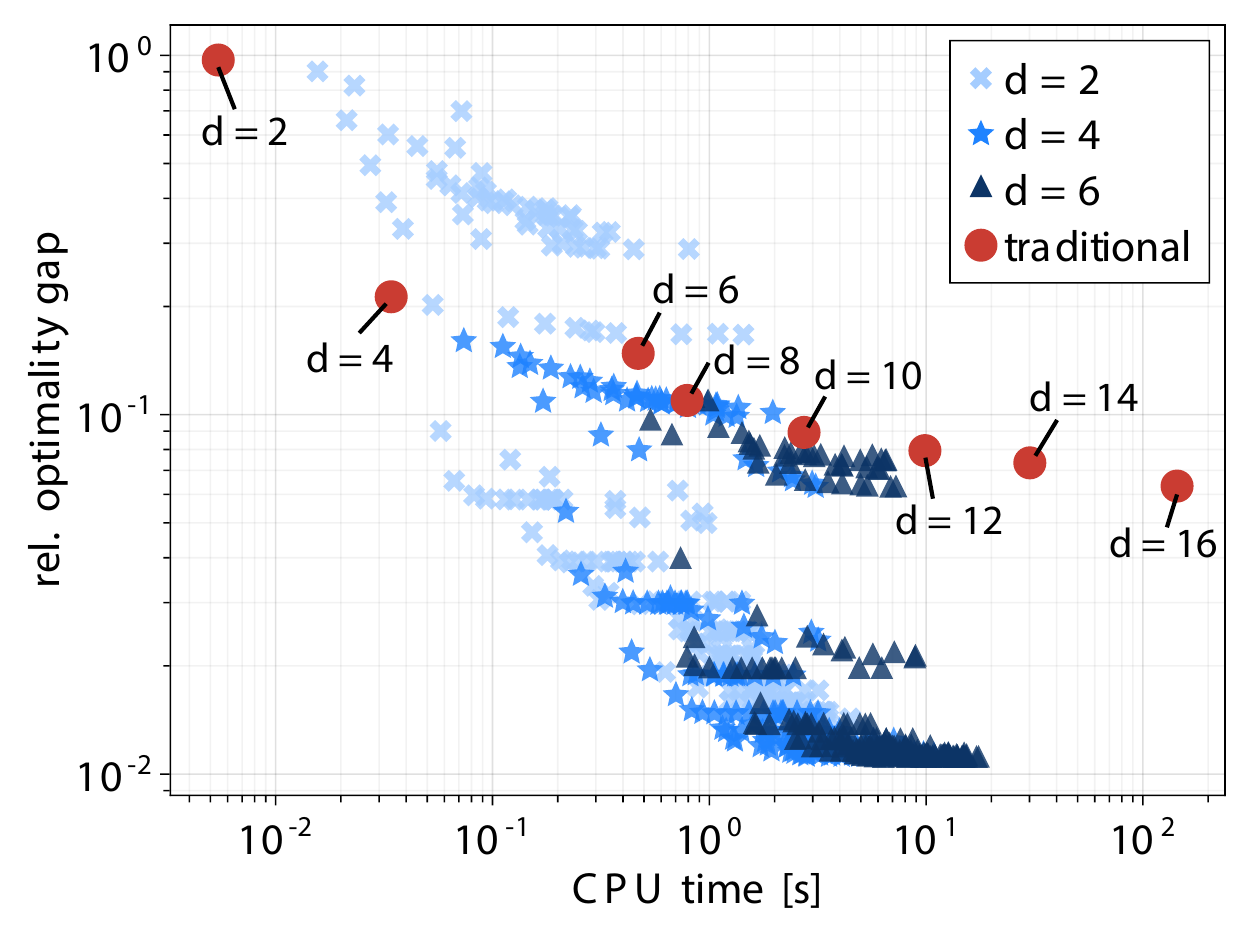}
    \caption{Trade-off between computational cost and bound quality for different approximation orders $d$ and domain partitions (different markers).}
    \label{fig:biocircuit_results}
\end{figure}

Figure \ref{fig:biocircuit_results} shows the trade-off between computational cost and bound quality achieved by different choices for the partition of the problem domain and approximation order. The bound quality is again measured by the relative optimality gap, estimated as described in Section \ref{sec:boundquality}. 
Analogous to the diffusion control example considered in Section \ref{sec:case_study}, the results demonstrate that an adequate partitioning of the problem domain substantially reduces the cost of computing bounds of a given quality when compared to the traditional approach. Moreover, notably tighter bounds could be computed overall due to a less conservative overapproximation of the process' infinite state space in the formulation of the bounding problems. 


\section{Conclusion}\label{sec:conclusion}
We have proposed a simple partitioning strategy for improving the practicality of MSOS relaxations for stochastic optimal control problems with polynomial data. From the primal perspective, this strategy can be interpreted as constructing the MSOS relaxation for an infinite-dimensional LP over finitely many occupation measures ``localized'' on elements of a partition of the control problem's space-time domain. From the dual perspective, the bounding problems seek a maximal piecewise-polynomial underestimator to the value function via sum-of-squares programming.

The key advantage of this framework over application of the MSOS hierarchy to the traditional occupation measure formulation for stochastic optimal control is that it offers a flexible and interpretable mechanism to tighten the obtained semidefinite bounding problems without degree augmentation -- simple refinement of the problem domain partition. On the one hand, this enables tightening of the bounding problems at as benign as linearly increasing cost, contrasting the combinatorial scaling incurred by naive degree augmentation. On the other hand, it promotes practicality by providing a way to avoid high degree sum-of-squares constraints and their notorious implications for poor numerical conditioning. As demonstrated with two examples, these advantages can lead to notable improvements in practical utility of the occupation measure approach to stochastic optimal control. 

In future work, we will investigate the use of distributed optimization techniques to further improve efficacy of the proposed framework by exploiting the weakly-coupled block structure of the bounding problems.

\section*{Acknowledgements}
FH gratefully acknowledges helpful discussions with Paul I. Barton and Frank Sch{\"a}fer. 

This material is based upon work supported by the National Science Foundation
under this material is based upon work supported by the National Science Foundation under grant no. OAC-1835443, grant no. SII-2029670, grant no. ECCS-2029670, grant no. OAC-2103804, and grant no. PHY-2021825. We also gratefully acknowledge financial support from the MIT-Switzerland Lockheed Martin Seed Fund and MIT International Science and Technology Initiatives (MISTI).

\section*{References}
\printbibliography

\appendix 
\section{Proof of Corollary 1} \label{app:cor1}

\begin{proof}
    Fix $z \in X$ and $t \in [t_{n_T-1},T)$. Now consider an admissible control process $u_s$ such that all paths of the controlled process $(s, x_s,u_s)$ lie in $[t,T] \times X \times U$ with $x_t \sim \delta_z$. Further define $\tau_0 = t$ and $\tau_i$ for $i \geq 1$ to be the minimum between $T$ and the time point at which the process crosses for the $i$\textsuperscript{th} time from one subdomain of the partition $X_1,\dots, X_{n_X}$ to another. 
    By construction, the process is confined to some (random) subdomain $X_k$ in the interval $[\tau_i, \tau_{i+1}]$. Since $w_{n_T,k}$ is sufficiently smooth on $[\tau_i, \tau_{i+1}] \times X_k$, Ito's lemma applies and yields that
    \begin{align*}
        w_{n_T,k}(\tau_{i+1}, x_{\tau_{i+1}}) = w_{n_T,k}(\tau_{i}, x_{\tau_i}) + \int_{\tau_i}^{\tau_{i+1}} \A w_{n_T, k} (s, x_s,u_s) \, {\rm d}s
        + \int_{\tau_i}^{\tau_{i+1}} \nabla_x w_{n_T,k}(s,x_s)^\top g(x_s,u_s) \, {\rm d} B_s.
    \end{align*}
    Now note that by Constraint \eqref{eq:path}, 
    \begin{align*}
        \int_{\tau_i}^{\tau_{i+1}} \A w_{n_T, k} (s, x_s,u_s) \, {\rm d}s \geq  -\int_{\tau_i}^{\tau_{i+1}} \ell(x_s,u_s) \, {\rm d}s,
    \end{align*}
    and further
    \begin{align*}
        \E{\delta_z}{\int_{\tau_i}^{\tau_{i+1}} \nabla_x w_{n_T,k}(s,x_s)^\top g(x_s,u_s) \, {\rm d} B_s} = 0
    \end{align*}
    as the integrand is square-integrable by Assumption \ref{ass:finite_moments} and $\tau_i \leq \tau_{i+1}$ are stopping times with respect to the natural filtration \cite[Chapter 2, Proposition 1.1]{Ikeda2014sdes}. Thus, after taking expectations, we obtain
    \begin{align*}
        \E{\delta_z}{ w_{n_T,k}(\tau_{i}, x_{\tau_{i}}) }  \leq \E{\delta_z}{\int_{\tau_i}^{\tau_{i+1}} \ell(x_s,u_s) \, {\rm d}s + w_{n_T,k}(\tau_{i+1}, x_{\tau_{i+1}}) }.
    \end{align*}
    
    Moreover, continuity holds at any crossing between any distinct subdomains $X_k$ and $X_j$ due to Constraint \eqref{eq:boundary} such that
    \begin{align*}
       \E{\delta_z}{w(\tau_{i}, x_{\tau_i})}  = \E{\delta_z}{ w_{n_T,k}(\tau_{i}, x_{\tau_i}) } = \E{\delta_z}{w_{n_T,j}(\tau_i, x_{\tau_i})},
    \end{align*}
    when the process crosses from $X_k$ to $X_j$ at $\tau_i$. Now using that $\E{\delta_z}{w(\tau_0, x_{\tau_0})} = w(t,z)$, we obtain by summing over the time intervals $[\tau_0, \tau_1], \dots, [\tau_{N},\tau_{N+1}]$ that 
    \begin{align*}
        w(t,z) \leq  
        \E{\delta_z}{\int_{t}^{\tau_{N+1}} \ell(x_s,u_s) \, {\rm d}s + w(\tau_{N+1}, x_{\tau_{N+1}}) \, {\rm d}s }.
    \end{align*}
    Letting $N \to \infty$, it follows that 
    \begin{align*}
        w(t,z) \leq \E{\delta_z}{\int_{t}^{T} \ell(x_s,u_s) \, {\rm d}s + w(T, x_T)  }
    \end{align*}
    as $\tau_N \to T$ almost surely. Finally using that $w(T,x) \leq \phi(x)$ on $X$ due to Constraint \eqref{eq:transversality} and the fact that all results hold for any admissible control policy, we obtain the desired result $w(t,z) \leq V(t,z)$. 
    
    It remains to show that $w$ preserves the lower bounding property across the boundaries introduced by discretization of the time domain. To that end, note that by an analogous argument as before, we have for any $t \in [t_{i-1}, t_{i})$ that 
    \begin{align*}
          w(t,z) \leq \E{\delta_z}{\int_{t}^{t_i} \ell(x_s,u_s) \, {\rm d}s + \lim_{s \nearrow t_{i}} w(s, x_{s}) }. 
    \end{align*}
    Since Constraint \eqref{eq:time} implies that $\lim_{s \nearrow t_{i}} w(s, x) \leq w(t_{i}, x)$ on $X$, it finally follows by induction that $w(t,z) \leq V(t,z)$ for any $t \in [0,T]$ and $z \in X$.
\end{proof}

\section{Measure transport equation}\label{app:slack}

For the sake of brevity, we derive expressions for the slack measures in \eqref{eq:disc_weakOCP} here only for the special case of a bipartition $X_1$ and $X_2 = X \setminus X_1$ of the state space $X$ and leave the temporal domain unpartitioned. Accordingly, we only need to consider crossings of the boundary $\partial X_{12} = \partial X_1 \cap \partial X_2$ between the two partition elements. Generalizations of this derivation to more complicated spatio-temporal partitions are straightforward.

For a fixed, admissible control policy $u_t$, we let $\tau_n^+$ ($\tau_n^-$) be the $n$\textsuperscript{th} time of entrance (exit) of the controlled process $(t, x_t)$ to (from) $[0,T]\times X_1$. We denote the associated entrance (exit) measures with $\mu_n^+$ ($\mu_n^-$) and further decompose them according to
\begin{align*}
    &\mu_n^+ = \nu_n^+ + \pi_n^+ \text{ and }\mu_n^- = \nu_n^- + \pi_n^-.
\end{align*}
Here, the measures $\pi_n^+ \in \mathcal{M}_+((0,T] \times \partial X_{12})$ and $\pi_n^- \in \mathcal{M}_+([0,T) \times \partial X_{12})$ describe entrances and exits of the process across the spatial boundary $\partial X_{12}$, while the measures $\nu_n^+ \in \mathcal{M}_+(\lbrace 0 \rbrace \times X_1)$ and $\nu_n^- \in \mathcal{M}_+(\lbrace T \rbrace \times X_1)$ complementarily cover entrances and exits at the initial and terminal time, respectively. The associated duality brackets read
\begin{align*}
    &\langle w, \nu_n^+ \rangle = \E{\nu_0}{w(\tau_n^+, x_{\tau_n^+}) \mathds{1}_{\lbrace 0 \rbrace \times X_1}(\tau_n^+, x_{\tau_n^+})}\\
    &\langle w, \pi_n^+ \rangle = \E{\nu_0}{w(\tau_n^+, x_{\tau_n^+}) \mathds{1}_{(0,T] \times \partial X_{12}}(\tau_n^+, x_{\tau_n^+})}\\
    &\langle w, \nu_n^- \rangle = \E{\nu_0}{w(\tau_n^-, x_{\tau_n^-}) \mathds{1}_{\lbrace T \rbrace \times X_1}(\tau_n^-, x_{\tau_n^-})}\\
    &\langle w, \pi_n^- \rangle = \E{\nu_0}{w(\tau_n^-, x_{\tau_n^-}) \mathds{1}_{[0,T) \times \partial X_{12}}(\tau_n^-, x_{\tau_n^-})}
\end{align*}

Application of Dynkin's formula to the stopping times $\tau_n^+$ and $\tau_n^-$ yields for an observable $w \in C^{1,2}([0,T]\times X_1)$ that
\begin{align*}
    \langle w, \nu_n^- + \pi_n^- \rangle - \langle w, \nu_n^+ + \pi_n^+ \rangle = \mathbb{E}_{\nu_0}\left[ \int_{\tau^+_i}^{\tau_n^-} \mathcal{A}w(s,x_s,u_s) \,{\rm d}s \right]
\end{align*}
and thus it follows for any $N \in \mathbb{N}$ that
\begin{align*}
    \langle w, \sum_{i=1}^N \nu_n^- \rangle - \langle w, \sum_{i=1}^N \pi_n^+ - \pi_n^- \rangle - \langle w, \sum_{i=1}^N \nu_1^+ \rangle = \mathbb{E}_{\nu_0}\left[ \sum_{i=1}^N \int_{\tau^+_i}^{\tau_n^-} \mathcal{A}w(s,x_s,u_s) \,{\rm d}s \right].
\end{align*}%
In the limit $N \to \infty$, the right-hand side of the above expression converges to the integral of $\mathcal{A}w$ with respect to the expected state-action occupation measure localized at $X_1$. Moreover, by construction $\nu_1^+$ coincides with the initial instantenous occupation measure localized on $X_1$, whereas $\nu_n^+ = 0$ holds for all $i > 1$. Similarly, the terminal local instantaneous occupation measure coincides with $\sum_{n=1}^{\infty} \nu_n^-$ as the process is by Assumption \ref{ass:finite_moments} non-explosive on $[0,T]$. Finally, it remains to observe that by construction
\begin{align*}
    \langle w, \sum_{n=1}^\infty \pi_n^+ - \pi_n^- \rangle =\E{\nu_0}{\sum_{n=1}^{N^{12}_+} w(\tau_{n+}^{12}, x_{\tau_{n+}^{12}}) - \sum_{n=1}^{N^{12}_-}  w(\tau_{n-}^{12}, x_{\tau_{n-}^{12}}),} 
\end{align*}
where $N^{12}_+$ ($N^{12}_-$) denote the (random) number of crossings of $\partial X_{12}$ into (out of) $X_1$ in $(0,T)$ at the corresponding times $\tau_{n+}^{12}$ ($\tau_{n-}^{12}$).

\end{document}